\documentclass[11pt]{article}

\usepackage{tikz}
\usepackage{subfigure}
\usepackage[english]{babel}

\usepackage[center]{caption2}
\usepackage{amsfonts,amssymb,amsmath,latexsym,amsthm}
\usepackage{multirow}
\usepackage[usenames,dvipsnames]{pstricks}
\usepackage{epsfig}
\usepackage{pst-grad} 
\usepackage{pst-plot} 
\usepackage[space]{grffile} 
\usepackage{etoolbox} 
\makeatletter 
\patchcmd\Gread@eps{\@inputcheck#1 }{\@inputcheck"#1"\relax}{}{}

\def\qedB{{\hfill\enspace\vrule height8pt depth0pt width8pt}}

\def\ex{\mathrm{ex}}

\begin{document}

\title{\bf\Large On codegree Tur\'an density of the 3-uniform tight cycle $C_{11}$}

\date{}

\author{
Jie Ma\footnote{School of Mathematical Sciences, University of Science and Technology of China, Hefei, Anhui 230026, China.
Research supported by National Key Research and Development Program of China 2023YFA1010201 and National Natural Science Foundation of China grant 12125106.
Email: jiema@ustc.edu.cn.}
}

\maketitle

\begin{abstract}
Piga, Sanhueza-Matamala, and Schacht recently established that the codegree Tur\'an density of 3-uniform tight cycles $C_\ell$ is $1/3$ for $\ell\in \{10, 13, 16\}$ and for all $\ell\geq 19$.
In this note, we extend their proof to determine the codegree Tur\'an density of the 3-uniform tight cycle $C_{11}$, thereby completing the picture for tight cycles of length at least 10.
\end{abstract}



Let $H$ be a 3-uniform hypergraph.
For any pair $\{u,v\}\subseteq V(H)$, the {\it codegree} of this pair, denoted by $d_H(uv)$, is the number of vertices $w\in V(H)$ satisfying $uvw\in E(H)$.
Let $\delta_2(H)=\min_{\{u,v\}\subseteq V(H)} d_H(uv)$ be the {\it minimum codegree} of $H$.
The {\it codegree Tur\'an number} $\ex_2(n,F)$ of a 3-uniform hypergraph $F$ denotes the maximum value of $\delta_2(H)$ among all $n$-vertex 3-uniform hypergraphs $H$ that do not contain a copy of $F$ as a subgraph.
The {\it codegree Tur\'an density} of $F$ is then defined to be the limit $\gamma(F)=\lim_{n\to \infty} \frac{\ex_2(n,F)}{n}.$
This notion was introduced by Mubayi and Zhao \cite{MZ} and has since attracted considerable attention in recent research (for detailed discussions, see \cite{FPVV, R20}).
Notably, Falgas-Ravry, Pikhurko, Vaughan, and Volec \cite{FPVV} proved that $\gamma(K_4^-)=1/4$, where $K_4^-$ denotes the hypergraph obtained from the complete 3-uniform hypergraph on four vertices $K_4$ by deleting an edge.
A conjecture of Czygrinow and Nagle \cite{CN} states that $\gamma(K_4)=1/2$, which remains open.

For integers $\ell\geq 5$, the {\it tight cycle} $C_\ell$ is defined as a 3-uniform hypergraph with vertex set $\{v_1,v_2,...,v_\ell\}$ and edge set $\{v_{i-1}v_iv_{i+1}: i \mbox{ is taken modulo } \ell\}$.
If $\ell$ is divisible by 3, then it is known that $\gamma(C_\ell)=0$.
One of the results in Balogh, Clemen, and Lidick\'y \cite{BCL} demonstrates that $\gamma(C_\ell)\leq 0.3993$ for all $\ell\ge 5$ except $\ell=7$.
On the other hand, there are constructions showing $\gamma(C_\ell)\geq 1/3$ for all $\ell$ not divisible by 3 (see, e.g., \cite{PSS}).
Very recently, using an elegant short proof, Piga, Sanhueza-Matamala, and Schacht \cite{PSS} determined the precise value of $\gamma(C_\ell)$ for all but finitely many choices of $\ell$.

\medskip

{\bf \noindent Theorem 1.} (Piga, Sanhueza-Matamala, and Schacht \cite{PSS}) For $\ell\in \{10, 13, 16\}$ and for every $\ell\geq 19$ not divisible by 3, it holds that $\gamma(C_\ell)=1/3$.

\medskip

As demonstrated by the authors in \cite{PSS}, it holds that $\gamma(C_{\ell+3})\leq \gamma(C_\ell)$ for all $\ell\geq 4$ and $\gamma(C_{t\ell})\leq \gamma(C_\ell)$ for all $t\geq 2$.
Therefore, to prove the aforementioned theorem, the authors in \cite{PSS} only needed to show that $\gamma(C_{10})\leq 1/3$.
They also speculated that $\gamma(C_\ell)=1/3$ for every $\ell\geq 5$ not divisible by 3.
In this note, we extend the arguments from \cite{PSS} to establish that $\gamma(C_{11})\leq 1/3$.
The following is our main result.

\medskip

{\bf \noindent Theorem 2.} $\gamma(C_{11})=1/3$, and consequently, $\gamma(C_{14})=\gamma(C_{17})=1/3$.

\medskip

Combining this with Theorem~1, $\gamma(C_\ell)=1/3$ holds for every $\ell\geq 10$ not divisible by 3.

We now give the proof of Theorem 2.
For two hypergraphs $G$ and $H$,
a {\it homomorphism} from $G$ to $H$ is a mapping $f:V(G)\to V(H)$ such that $f(e)\in E(H)$ for every $e\in E(G)$.

\medskip

{\noindent \bf Proof of Theorem 2.}
This proof builds on the approach in \cite{PSS}, with the addition of one key idea: to analyze the types of the three pairs in each edge.

Consider an arbitrary $\epsilon>0$ and sufficiently large integers $n$.
Let $H$ be any 3-uniform hypergraph on $n$ vertices with $\delta_2(H)\geq (1/3+\epsilon)n$.
It suffices to prove that $H$ contains an homomorphic copy of the tight cycle $C_{11}$.
Following the notion given in \cite{PSS}, we call the only vertex with degree 3 in a $K_4^-$ as the {\it apex} of that $K_4^-$.
A pair $\{u,v\}$ of distinct vertices in $H$ is called an {\it apex pair} if there exists a $K_4^-$ of $H$ containing $u$ and $v$ with the apex being either $u$ or $v$.
If $\{u,v\}$ is an apex pair with apex $v$, then we denote this relationship by the arc $u\to v$.
Let $D$ be the digraph consisting of all such arcs $u\to v$.
A pair $\{u,v\}$ is a {\it base pair} if there exists a $K_4^-$ of $H$ containing $u$ and $v$ such that neither $u$ nor $v$ is the apex.

First we claim that every edge $xyz$ in $H$ is contained in a $K_4^-$ of $H$, and moreover, $D[\{x,y,z\}]$ has a vertex with indegree 2.
Since $\delta_2(H)\geq (1/3+\epsilon)n$, we see that $d_H(xy)+d_H(xz)+d_H(yz)\geq (1+3\epsilon)n$.
So there exists a vertex $w$ belonging to at least two of the neighbourhoods $N_H(xy), N_H(xz)$ and $N_H(yz)$ (say the former two).
Then $\{xyz,xyw,xzw\}$ induces a $K_4^-$ of $H$ with the arcs $y\to x$ and $z\to x$, proving the claim.

Next, we show that if a pair $\{x,y\}$ is both an apex pair and a base pair, then $H$ contains an homomorphic copy of $C_{11}$.
Let $K$ be a $K_4^-$ of $H$ with $V(K)=\{x,y,a,b\}$ and apex $x$.
Additionally, let $K'$ be another $K_4^-$ of $H$ with $V(K')=\{x,y,c,d\}$ and apex $c$.
We observe that the sequence $(x,c,d,y,c,x,y,b,x,a,y)$ forms an homomorphic copy of $C_{11}$.

We also claim that if $xyz$ is an edge in $H$ with arcs $y\to x$ and $z\to y$, then $H$ contains an homomorphic copy of $C_{11}$.
To see this, let $K$ be a $K_4^-$ of $H$ with $V(K)=\{x,y,a,b\}$ and apex $x$,
and let $K'$ be a $K_4^-$ of $H$ with $V(K')=\{y,z,c,d\}$ and apex $y$.
Then the sequence $(x,a,b,x,y,z,d,y,c,z,y)$ forms an homomorphic copy of $C_{11}$ in $H$.

Recall the digraph $D$.
Define $B$ to be the graph consisting of all base pairs.
Together with the previous statements, we can conclude that
\begin{itemize}
\item [(1).] There are no 2-cycles in $D$, and
\item [(2).] For every edge $xyz$ in $H$, there exists a vertex, say $x$, such that
$D[\{x,y,z\}]$ has exactly two arcs $y\to x$ and $z\to x$, and $B[\{x,y,z\}]$ has a unique edge $yz$.
\end{itemize}
Using these two items, we can conclude the following (consistent with Claim 4 of \cite{PSS}):
\begin{itemize}
\item [(a).] If $d_B(v)>0$, then $d^+_D(v)\geq (1/3+\epsilon)n$.
\item [(b).] If $d_D^+(v)>0$, then $d_B(v)\geq (1/3+\epsilon)n$.
\item [(c).] If $d_D^-(v)>0$, then $d^-_D(v)\geq (1/3+\epsilon)n$.
\end{itemize}
To see item (a), let $uv\in E(B)$.
For any $w$ with $uvw\in E(H)$, the above item (2) implies that $v\to w$.
Thus $d^+_D(v)\geq \delta_2(H)\geq (1/3+\epsilon)n$.
The other items can be derived similarly.

Now we find ourselves in exactly the same situation as in the proof of \cite{PSS}.
Referring to the last two paragraphs of that proof, we can demonstrate that there exists a vertex $v^*$ with $d_B(v^*)>0$, $d_D^+(v^*)>0$, and $d_D^-(v^*)>0$.
Using items (1), (a), (b) and (c), we arrive at the final contradiction:
the neighborhoods $N_B(v^*), N_D^+(v^*), N_D^-(v^*)$ are disjoint, and their union exceeds the total number of vertices.
\qedB

\medskip

\bigskip

\noindent {\bf Acknowledgements.}
I would like to thank Xizhi Liu, Tianhen Wang, and Tianming Zhu for their helpful discussions and for carefully reading a draft.
I also express my gratitude for the hospitality extended to me during my visit to the School of Mathematics and Statistics at Shandong University, Weihai, in August 2024, which inspired part of this research.

\end{document}